\documentclass[leqno]{amsart}
\usepackage{amsmath, amssymb}
\usepackage{amsfonts}
\usepackage{graphics}
\usepackage[dvips]{graphicx}
\usepackage{eepic}
\usepackage{epic}
\newtheorem{thm}{Theorem}[section]
\newtheorem{crl}[thm]{Corollary}

\newtheorem{lmm}[thm]{Lemma}

\newtheorem{prp}[thm]{Proposition}

\newtheorem{rem}{Remark}

\begin{document}
\begin{large}
\title{Asymptotics of degenerating Eisenstein series }
\author{Kunio Obitsu\\
Faculty of Science, Kagoshima University}
\thanks{The author is partially supported by JSPS Grant-in-Aid for
Exploratory Research 2005-2007. {\it Mathematical Subject Classification (2000)}: 11M36, 32G15, 53C43.} 
\maketitle

\begin{center}
\thanks{Dedicated to Takahide Kurokawa and Kimio Miyajima\\
on the occasion of their 60th birthdays }     
\end{center}
   
\begin{abstract}      
We give some estimates for the asymptotic orders of degenerating Eisenstein series
for certain families of degenerating punctured Riemann surfaces, motivated
by the question of identifying $L_{2}$-cohomology of the Takhtajan-Zograf metric that
is originally asked by To and Weng. 
\end{abstract}

\section{ Introduction}

We consider the Teichm\"uller space $T_{g,n}$ and the associated Teichm\"uller 
curve ${\mathcal T}_{g,n}$ of Riemann surfaces 
of type $(g,n)$ (i.e., Riemann surfaces of genus $g$ and with $n>0$ 
punctures). We will assume that $2g-2+n>0$, so that 
each fiber of the holomorphic projection map $\pi:{\mathcal T}_{g,n}\to T_{g,n}$ is 
stable or equivalently, it admits the 
complete hyperbolic metric of constant sectional curvature $-1$. The kernel of 
the differential $T{\mathcal T}_{g,n}\to 
TT_{g,n}$ forms the so-called vertical tangent bundle over ${\mathcal T}_{g,n}$, 
which is denoted by $T^{V}{\mathcal T}_{g,n}$. 
The hyperbolic metrics on the fibers induce naturally a Hermitian metric on 
$T^{V}{\mathcal T}_{g,n}$. 

In the study of the family of $\bar\partial_k$-operators acting on the 
$k$-differentials on Riemann surfaces (i.e., 
cross-sections of $\big(T^{V}{\mathcal T}_{g,n}\big)^{-k}\big|_{\pi^{-1}(s)}\to 
\pi^{-1}(s),\ s\in T_{g,n}$), Takhtajan 
and Zograf introduced in  {\cite{TZ}} a K\"ahler metric on $T_{g,n}$, 
which is known as the Takhtajan-Zograf metric. In 
 {\cite{TZ}}, they showed that the Takhtajan-Zograf metric is invariant under the 
natural action of the Teichm\"uller modular 
group $\text{Mod}_{g,n}$ and it satisfies the following remarkable identity 
on $T_{g,n}$: 
\begin{equation*}
c_1(\lambda_k, \left\| \cdot \right\|_{Q,k})={{6k^2-6k+1}\over {12 \pi^2}}\
\omega_{\text{WP}}-{1\over 9}\omega_{\text{TZ}}.
\end{equation*}
Here $\lambda_k=\text{det}(\text{ind}\,\bar\partial_k)=\bigwedge^{\text {max}}
\text{Ker}\,\bar\partial_k\otimes 
(\bigwedge^{\text {max}}\text{Coker}\,\bar\partial_k)^{-1}$ denotes the determinant 
line bundle on $T_{g,n}$, $\left\| \cdot \right\|_{Q,k}$ 
denotes the Quillen metric on $\lambda_k$, and $\omega_{\text{WP}}$, 
$\omega_{\text{TZ}}$ denote the K\"ahler form of 
the Weil-Petersson metric, the Takhtajan-Zograf metric on $T_{g,n}$ 
respectively.  In  {\cite{We}}, Weng studied the Takhtajan-Zograf metric in terms of Arakelov intersection, and he proved that $\frac{4}{3} \omega_{\text{TZ}}$ coincides with the first Chern form of an associated metrized Takhtajan-Zograf line bundle over the moduli space $\mathcal M_{g,n}=T_{g,n}/\text{Mod}_{g,n}$.  Recently, Wolpert  {\cite{Wo3}} gave a natural definition 
of a Hermitian metric on the Takhtajan-Zograf line bundle whose first Chern form gives  
$\frac{4}{3} \omega_{\text{TZ}}$.
Furthermore, we can observe that in the second term of the asymptotic expansion of
the Weil-Petersson metric near the boundary of $\mathcal M_{g,n}$, the Takhtajan-Zograf
metrics on the boundary moduli spaces could appear (see {\cite{OW}}).

We propose a program of identifying 
$L_{2}$-cohomology of $\mathcal M_{g,n}$
with respect to the Takhtajan-Zograf metric $H^{*}( \mathcal M_{g,n}, \omega_{TZ})$. 
Originally, Saper ({\cite{S}}) applied Masur's formula
({\cite{M}}) to show that 
$L_{2}$-cohomology of $\overline{\mathcal M}_{g,0}$
with respect to the Weil-Petersson metric
$H^{*}( \mathcal M_{g,0}, \omega_{WP})$
 is naturally isomorphic to 
$H^{*}( \overline{\mathcal M}_{g,0}, {\mathbf R} )$. However, 
it is disappointing that the results for the asymptotics of the Takhtajan-Zograf
metrics in {\cite{OTW1}} are not sufficient for us to determine $H^{*}( \mathcal M_{g,n}, \omega_{TZ})$.

 In the present paper, we prove some estimates for the degenerating 
orders of Eisenstein series for certain families of degenerating punctured Riemann surfaces, which may be an important step for calculating $H^{*}( \mathcal M_{g,n}, \omega_{TZ})$. It should be noted that there are already some results 
for the behaviors of degenerating Eisenstein series 
({\cite{F}}, {\cite{G}}, {\cite{O}}, {\cite{Sc}}).

The author would like to thank K. Matsuzaki for showing him properties
of thick parts of Riemann surfaces.
He would like to thank S.A. Wolpert for showing him properties for 
the modified harmonic map. 
Furthermore, he would like to thank W.-K. To and L. Weng for
posing the problem to identify $L_{2}$-cohomology of $\mathcal M_{g,n}$
with respect to the Takhtajan-Zograf metric. He is grateful to the referee for his careful
reading.

\section{Main Theorems}

\subsection{Settings and notation}

 For simplicity of exposition, we consider a degenerating family $\{S_{l}\}$  of
Riemann surfaces of type $(g,1)$ with two zero-homologous 
pinching geodesics $\gamma_{1}$ and $\gamma_{2}$ which divide the surface $S_{l}$
into three components $S_{l}^{1}, S_{l}^{2}, S_{l}^{3}$:
the geodesic $\gamma_{1}$ divides $S_{l}^{1}$ from $S_{l}^{2}$, the geodesic  $\gamma_{2}$
divides $S_{l}^{2}$ from $S_{l}^{3}$,
and $S_{l}^{1}$ has the unique puncture.
(It should be noted that all claims in propositions, theorems, etc. are easily generalized to the case of  any degenerating family of hyperbolic surfaces of finite type with at least one puncture. In some of the statements, we will give remarks for the general case.)
The vector-valued parameter $l$ varies around the origin in the Euclidian space 
$\mathbf{R}^{6g-4}$,
where $l=0$ represents the unique degenerate surface $S_{0}$ in the family.
The limit surface $S_{0}$ consists of three components 
$S_{0}^{1}$, $S_{0}^{2}$, $S_{0}^{3}$ which are the limits of 
$S_{l}^{1}, S_{l}^{2}, S_{l}^{3}$ respectively as $l \rightarrow 0$. 
Let $q_{j}$ be the node shared
by $S_{0}^{j}$ and $S_{0}^{j+1}$ ($j=1,2$).
A puncture the smooth surface in the degenerate family originally has 
will be called an $\it{old}$ puncture, for simplicity.
It should be noted that the degenerating family can be described by the modified 
infinite-energy harmonic maps 
$f^{l}: S_{0} \longrightarrow S_{l}\backslash \{ \gamma_{1}, \gamma_{2} \}$, which
are introduced by S. Wolpert (\cite{Wo2}).

Let $L_{l}(\gamma)$ be the hyperbolic length of a simple closed geodesic $\gamma$ on $S_{l}$.
For $0 \leq k \le 1$ and $j=1,2$, set
$$N_{\gamma_{j}}(k) = \bigg\{  p \in S_{l} \biggm| d_{l}(p, \gamma_{j}) 
\leq k\ \text{sinh}^{-1} \bigg(1 \Big/
\text{sinh}\frac{ L_{l}(\gamma_{j})}{2}\bigg) \bigg\},$$
 the collar neighborhood around $\gamma_{j}$ in $S_{l}$, where $d_{l}(\cdot, \cdot)$ denotes the hyperbolic distance on $S_{l}$. Here we remark that
\begin{equation*}(2.1)   \hspace{2.0cm}
\text{sinh}^{-1} \Big(1 \big/\text{sinh}\frac{x}{2} \Big)
= - \log x +2\log 2 + O(x^{2}), \quad x \rightarrow 0,
\end{equation*}which will be essentially used in the proofs of Lemma 2.3 and Lemma 2.6.

For $a \ge 1$, the $a$-cusp region $ C_{j}(a)$ ($\subset  S_{0}^{j} \cup S_{0}^{j+1}$)
around the node $q_{j}$  is the union of two 
copies of $\langle z \mapsto z+1 \rangle
\backslash \{z \in H \;|\; \operatorname{Im}z \ge a \} $, equipped with the metric $ds^{2} = (dy^{2} +
dx^{2}) / y^{2},$ where $H:= \{z \in \mathbf{C}\; |\; \operatorname{Im} z >0 \}$ is 
the upper half plane.

Let $(f^{l})^{*} \Delta_{l}$ denote the pull-back of the negative hyperbolic 
Laplacian $\Delta_{l}$ on $S_{l}$ by $f^{l}$, that is, for a $C^{2}$-function $h$ on $S_{0}$,
$$ (f^{l})^{*} \Delta_{l}\; (h) = \Delta_{l} \bigm( h \circ (f^{l})^{-1} \bigm) \circ f^{l}.$$
Let $\Delta_{0}$ denote the negative hyperbolic Laplacian on $S_{0}$. Then, it is known that
  $(f^{l})^{*} \Delta_{l}$ converges to $\Delta_{0}$
uniformly on any compact subset of $S_{0}$ in the $C^{3}$-norm (see  {\cite{Wo2}} ).
And, for a function  $g$ on $S_{l}$, the pull-back of $g$ by $f^{l}$
is defined as
$$ (f^{l})^{*}g = g \circ f^{l}.$$
It should be noted that a $C^{2}$-function $g$ on $S_{l}$ satisfies
\begin{align*} \hspace{-3.0cm}  (2.2)  \hspace{3.0cm}  
(f^{l})^{*} \Delta_{l} \bigm( (f^{l})^{*}g \bigm)  = \Delta_{l} (g) \circ f^{l},
\end{align*}
which will be used in the proof of Lemma 2.7.
 
\subsection{The counting function of orbits}

 Let $\Gamma_{l}$ be a Fuchsian group uniformizing $S_{l}$ such that 
$S_{l} \simeq H/ \Gamma_{l}$.
We normalize it such that $\Gamma_{l}$ contains a parabolic element 
$z \mapsto z+1$. A cyclic group generated by the parabolic element is denoted by
 $\Gamma_{\infty}$. 
Then the Eisenstein series for $\Gamma_{l}$ associated to the unique puncture is expressed as

$$E^{l}(z,s) = \sum\limits_{\delta \in \Gamma_{\infty} \backslash \Gamma_{l}} 
(\operatorname{Im} \delta z)^{s}, \quad z \in H,\;
\operatorname{Re} s >1.$$
Here for any $z$ in $H$ and any equivalent class $[\delta]$ in $\Gamma_{\infty} \backslash \Gamma_{l}$, we can select the unique representative $\hat{\delta}$ for $[\delta]$ such that $- \frac{1}{2} \leq \operatorname{Re} \hat{\delta} z < \frac{1}{2}$. Such $\hat{\delta}=\hat{\delta} (z, [\delta])$ 
will be called the {\it canonical} representative. 

$E^{l}(z,s) $ is invariant under the action of $\Gamma_{l}$. Thus it can be considered
as a function on $S_{l}$. Moreover, it is well known that the Eisenstein seires satisfies 
\begin{align*} \hspace{-1.0cm}(2.3) \hspace{2.0cm}
(\Delta - s(s-1))\;  E^{l}(z,s) = 0 , \quad  z \in H,\; 
\operatorname{Re} s >1, 
\end{align*}
which will play a crucial role in the proof of Lemma 2.7.
Here $\Delta := 4\; (\operatorname{Im} z)^{2}
 \frac{\partial^{2}}{\partial z \partial \bar{z}}$ is the negative hyperbolic Laplacian on $H$ 
 , invariant under $\Gamma_{l}$, and thus it naturally descends to $\Delta_{l}$
 on $S_{l}$. 

Now we are ready to present a new way to study the asymptotics of the Eisenstein series.
When $\operatorname{Im} z < 1$ and
$z$ is not equivalent to any point of $\{ w \in H |\ \operatorname{Im} w >1 \}$
under the action of $\Gamma_{\infty} \backslash  \Gamma_{l}$, it is easy to see that 
for $[\delta]$ in $\Gamma_{\infty} \backslash  \Gamma_{l}$, $\operatorname{Im} \hat{\delta}(z)= 
e^{-d (h, \hat{\delta} z)},$
where $d(\cdot, \cdot)$ denotes the hyperbolic distance in $H$ and 
$h= \{w \in H | - \frac{1}{2} \leq \operatorname{Re} w < \frac{1}{2}, \operatorname{Im} w =1 \}$.

We introduce two counting functions of orbits of $z$ with 
 $\operatorname{Im} z < 1$,
\begin{align*}
\Pi_{l} (h,z,t) & := \sharp \{ [\delta] \in \Gamma_{\infty} \backslash \Gamma_{l} \;|\
d(h, \hat{\delta} z) \leq t \},\\
\Pi_{l} (z,t) & := \sharp \{ [\delta] \in \Gamma_{\infty} \backslash \Gamma_{l} \;|\
d(i, \hat{\delta} z) \leq t \},
\end{align*}
where $\hat{\delta}$ is the {\it canonical} representative. 
Here we should remark that $d(i, \hat{\delta} z) =
 \min\limits_{\delta \in [\delta]} d(i, \delta z) $.

For $z$ with $\operatorname{Im} z < 1$ not equivalent to any point of 
$\{ w \in H |\ \operatorname{Im} w >1 \}$
under the action of $\Gamma_{\infty} \backslash \Gamma_{l}$, we can observe 
$$E^{l}(z,s) = \int_{0}^{\infty} e^{-st}\ d \Pi_{l} (h,z,t).$$

We will state a famous property of $\Pi_{l} (z,t)$ as in the form suited to our 
purpose.

\begin{prp}
There exists an absolute constant $U$ such that for $z \in H$ with 
$ \operatorname{Im}z <1$,
the following estimate holds:
$$\Pi_{l}(z,t) \leq U e^{t} \qquad \text{for any}\ t \geq 0
\; \text{and any}\ \Gamma_{l}.$$
\end{prp}
\begin{proof}
Our proof is based on the discussion in \cite{T} p.516.
Let $B(p,r)$ denote a hyperbolic ball centered at $p$ with radius $r$ in $H$.
Now the collar lemma assures us that in any hyperbolic surface with at least one puncture, each 
puncture has  a horocyclic neighborhood with area $2$  (see {\cite{SS}}).
Then we can find a universal constant $\varepsilon > 0$ such 
that orbits $B(\delta i, \varepsilon)$ for $\delta \in \Gamma_{l}$
are mutually disjoint for any $\Gamma_{l}$.
Because if $d(\delta i, z) \leq t$ then $B(\delta i,\varepsilon) \subset B(z,t+\varepsilon)$,
we have
\begin{align*}
\Pi_{l}(z,t) &= 
\sharp \{ [\delta] \in \Gamma_{\infty} \backslash \Gamma_{l} \;|\
d( (\hat{\delta})^{-1} i, z) \leq t \}\\
& \leq 
\sharp \{ [\delta] \in \Gamma_{\infty} \backslash \Gamma_{l} \;|\
(\hat{\delta})^{-1} (B(i, \varepsilon) )\subset B(z, t+\varepsilon)  \}\\
& \leq \frac{|B(z, t+\varepsilon)|}{|B(i, \varepsilon)|} =
\sinh^{2} \Bigl( {\frac{t+\varepsilon}{2}} \Bigr) \Bigl / \sinh^{2}{\frac{\varepsilon}{2}}\\
& \leq \frac{e^{\varepsilon}}{2\sinh^{2}{\frac{\varepsilon}{2}}}e^{t}
 \quad \text{for}\ t \geq 0.
\end{align*}
Here $| \cdot|$ denotes the hyperbolic area in $H$.
\end{proof}

\begin{prp}
Let $s >1$.
Let $z \in H$ with $ \operatorname{Im} z <1$ be not equivalent to any point 
of $\{ w \in H |\ \operatorname{Im}  w >1 \}$
under the action of $\Gamma_{\infty} \backslash \Gamma_{l}$.
Then we obtain
$$\Pi_{l} (z,t) \leq \Pi_{l} (h,z,t) \leq \Pi_{l} (z,t+1),$$
$$E^{l} (z,s) = s \int_{0}^{\infty} e^{-st} \Pi_{l} (h,z,t)\ dt,$$
$$s \int_{0}^{\infty} e^{-st}\ \Pi_{l} (z,t)\ dt \leq
E^{l} (z,s) \leq
s \int_{0}^{\infty} e^{-st}\ \Pi_{l} (z,t+1)\ dt.$$

\end{prp}
\begin{proof}
Because $d(i, \delta z) \leq d(h, \delta z) +1$, it follows from Proposition 2.1 that
$$\Pi_{l} (h,z,t) \leq \Pi_{l} (z,t+1) \leq e U e^{t}. $$
Then, integrations by parts and Proposition 2.1 provide
\begin{align*}
E^{l} (z,s) &=  \int_{0}^{\infty} e^{-st}\ d \Pi_{l} (h,z,t)
 = [e^{-st} \Pi_{l} (h,z,t)]_{0}^{\infty} 
 + s \int_{0}^{\infty} e^{-st} \Pi_{l} (h,z,t)\ dt\\
 & = s \int_{0}^{\infty} e^{-st} 
 \Pi_{l} (h,z,t)\ dt
 \leq s\int_{0}^{\infty} e^{-st} \Pi_{l} (z,t+1)\ dt.
\end{align*}
This is the right-hand inequality in the statement.

Next we will prove the left-hand inequality.
Since $d(h, \delta z) \leq d(i, \delta z)$, it is easy to see that
\begin{align*}
\Pi_{l} (h,z,t) &\geq \Pi_{l} (z,t),\\
E^{l} (z,s) &= 
s \int_{0}^{\infty} e^{-st} \Pi_{l} (h,z,t)\ dt\\
& \geq 
s \int_{0}^{\infty} e^{-st} 
 \Pi_{l} (z,t)\ dt.
\end{align*}
\end{proof}

\subsection{Upper bounds for degenerating Eisenstein series}

We are going to present upper bounds for Eisenstein series on the components 
$S_{l}^{2}$ and $S_{l}^{3}$.
\begin{lmm}
Assume $\operatorname{Re} s > 1$.
There exists an absolute constant $M_{2}(\operatorname{Re} s)$
depending only on 
 $\operatorname{Re} s$ such that 
 for $L_{l}(\gamma_{1}), L_{l}(\gamma_{2}) <
 2 \operatorname{sinh}^{-1}1$ and $0 \leq k \leq 1$, then
\begin{align*} 
& | E^{l}(z,s) | \le M_{2}(\operatorname{Re} s) \  
L_{l}(\gamma_{1})^{(1+k)(\operatorname{Re} s -1)} \qquad \qquad 
\qquad  \qquad \hspace{1.5pt} \text{on} \ \partial N_{\gamma_{1}}(k) \cap S_{l}^{2},\\
& | E^{l}(z,s) | \le M_{2}(\operatorname{Re} s)\ 
L_{l}(\gamma_{1})^{2(\operatorname{Re}s -1)}L_{l}(\gamma_{2})^{(1+k)(\operatorname{Re} s -1)} \qquad  
 \text{on} \ \partial N_{\gamma_{2}}(k) \cap S_{l}^{3}.
\end{align*}
\end{lmm}
\begin{proof}
Because $| E^{l}(z,s) | \leq E^{l}(z,\operatorname{Re} s)$ holds,
it is enough to show in the case $s > 1$.
 For $z$ in $H$,
$[z]$ denotes the corresponding point of $S_{l}$. 
$(2.1)$ implies easily that the distance of any curve connecting $[z]$ on 
$\partial N_{\gamma_{1}}(k) \cap S_{l}^{2}$  and the horocycle $[h]$ is greater than 
$(1+k) \text{[the width of half collar]}$. Therefore we see
$$\Pi_{l} (h,z,-(1+k) \log L_{l}(\gamma_{1})) = 0.$$
Then Proposition 2.2 yields
$$E^{l} (z,s) = s \int_{-(1+k) \log L_{l}(\gamma_{1})}^{\infty} e^{-st} \Pi_{l} (h,z,t)\ dt.$$
By Propositions 2.1 and 2.2, it concludes that
\begin{align*}
E^{l} (z,s) 
& \leq 
s \int_{-(1+k) \log L_{l}(\gamma_{1})}^{\infty} e^{-st} \Pi_{l} (z,t+1)\ dt\\
& \leq
s \int_{-(1+k) \log L_{l}(\gamma_{1})}^{\infty} e^{-st} eU e^{t}\ dt\\
& =
eUs \int_{-(1+k) \log L_{l}(\gamma_{1})}^{\infty} e^{-(s-1)t}\ dt\\
&=
\frac{eUs}{s-1}\ L_{l}(\gamma_{1})^{(1+k)(s-1)}.
\end{align*}
The second case is similar. Just replace $-(1+k) \log L_{l}(\gamma_{1})$ with
$-2 \log L_{l}(\gamma_{1})-(1+k) \log L_{l}(\gamma_{2})$.

\end{proof}

\begin{crl}
Assume as in Lemma 2.3.
Then for all $L_{l}(\gamma_{1}), L_{l}(\gamma_{2})< 2 \operatorname{sinh}^{-1}1 $ 
and all $k$ with $0 \leq k \leq 1$, it holds that
\begin{align*} 
& | E^{l}(z,s) | \le M_{2}(\operatorname{Re} s) \  L_{l}(\gamma_{1})^{(1+k)(\operatorname{Re} s -1)}
 \qquad \qquad \qquad \qquad \hspace{5.5pt}
 \text{on}\  S_{l}^{2}-N_{\gamma_{1}}(k),\\
& | E^{l}(z,s) | \le M_{2}(\operatorname{Re} s)\ 
L_{l}(\gamma_{1})^{2(\operatorname{Re}s -1)}\
L_{l}(\gamma_{2})^{(1+k)(\operatorname{Re} s -1)} \qquad  
 \text{on}\  S_{l}^{3}- N_{\gamma_{2}}(k).
\end{align*}
Here $M_{2}(\operatorname{Re} s)$ is the constant appearing in Lemma 2.3. 
\end{crl}
\begin{proof}
Because $| E^{l}(z,s) | \leq E^{l}(z,\operatorname{Re} s)$ holds,
it is enough to show the statements for $s > 1$.
By $(2.3)$, it is easy to see that $E^{l}(z,s)$ is subharmonic.
The maximal principle for subharmonic functions provides
\begin{align*}
\sup\limits_{z \in S^{2}_{l}- N_{\gamma_{1}}(k)} E^{l}(z,s) 
& \leq \sup\limits_{z \in S^{3}_{l} \cup S^{2}_{l}-N_{\gamma_{1}}(k)} E^{l}(z,s) \\
& =  \sup\limits_{z \in \partial N_{\gamma_{1}}(k) \cap S^{2}_{l}} E^{l}(z,s) \\
& \leq M_{2}(s)\ L_{l}(\gamma_{1})^{(1+k)(s -1)}.
\end{align*}
(Remark: even in the case $S^{3}_{l} \cup S^{2}_{l}$ has other {\it old} punctures, 
our discussion remains valid because $E^{l}(z,s)$ assumes $0$ at the {\it old} punctures.)
The second case is similar. Just use the second inequality in Lemma 2.3.
\end{proof}

We will summarize the special case for $k=0,1$ in Corollary 2.4 as follows.

\begin{thm}
Assume $\operatorname{Re} s > 1$.
Then for all $L_{l}(\gamma_{1}), L_{l}(\gamma_{2})< 2 \operatorname{sinh}^{-1}1 $, 
it holds that
\begin{align*} 
& | E^{l}(z,s) | \le M_{2}(\operatorname{Re} s) \  
L_{l}(\gamma_{1})^{(\operatorname{Re}s-1)} \qquad \qquad \qquad \hspace{30pt}
 \text{on}\  S_{l}^{2},\\
& | E^{l}(z,s) | \le M_{2}(\operatorname{Re} s) \  
L_{l}(\gamma_{1})^{2(\operatorname{Re}s-1)} \qquad \qquad \qquad \hspace{26pt}
 \text{on}\  S_{l}^{2} -N_{\gamma_{1}}(1),\\ 
& | E^{l}(z,s) | \le M_{2}(\operatorname{Re} s)\ 
L_{l}(\gamma_{1})^{2(\operatorname{Re}s -1)}\
L_{l}(\gamma_{2})^{(\operatorname{Re}s-1)} \qquad \hspace{4pt}
 \text{on}\  S_{l}^{3},\\
& | E^{l}(z,s) | \le M_{2}(\operatorname{Re} s)\ 
L_{l}(\gamma_{1})^{2(\operatorname{Re}s -1)}\
L_{l}(\gamma_{2})^{2(\operatorname{Re}s-1)} \qquad \text{on}\  
S_{l}^{3} - N_{\gamma_{2}}(1).
\end{align*}
Here $M_{2}(\operatorname{Re} s)$ is the constant appearing in Lemma 2.3. 
\end{thm}
\begin{rem}
Corollary 2.4 and Theorem 2.5 have essentially improved the order estimates for 
the degenerating Eisenstein series in {\cite{O}} Theorem 1 (2).
\end{rem}

\subsection{Lower bounds for degenerating Eisenstein series}

Now we are ready to present lower bounds for Eisenstein series on the components 
$S_{l}^{2}$ and $S_{l}^{3}$. Henceforth, the set of points in $S_{l}$ the injectivity radii of 
which are greater than $\sinh^{-1} 1$ will be called the {\it thick part} of $S_{l}$.

\begin{lmm}
Let $s > 1$. There exist positive constants $K_{i}=K_{i}(s, \{S_{l} \}) \; (i=1,2,3)$ 
depending only on $ s$ and the degenerating family $\{S_{l} \}$ such that 
 for $L_{l}(\gamma_{1}), L_{l}(\gamma_{2}) < 2 \operatorname{sinh}^{-1}1 $
 and $0 \leq k \leq 1$, then
 \begin{align*} 
&  E^{l}(z,s)  \ge K_{1} \  L_{l}(\gamma_{1})^{(1+k)s} \qquad \qquad 
\hspace{26pt}
 \text{on} \ \partial N_{\gamma_{1}}(k) \cap S_{l}^{2},\\
 &  E^{l}(z,s)  \ge 
 K_{2} \ L_{l}(\gamma_{1})^{2s} L_{l}(\gamma_{2})^{(1-k)s} \qquad \quad 
 \text{on} \ \partial N_{\gamma_{2}}(k) \cap S_{l}^{2},\\
&  E^{l}(z,s) \ge K_{3}\ L_{l}(\gamma_{1})^{2s}
   L_{l}(\gamma_{2})^{(1+k)s} \; \; \quad \quad  \hspace{5pt}
 \text{on} \ \partial N_{\gamma_{2}}(k) \cap S_{l}^{3}.
\end{align*}
\end{lmm}
\begin{proof}
We mimic the proof of Lemma 4.2 in \cite{Wo1} p.84.
For $z \in H$ with $ \operatorname{Im} z < 1$,
$$(\text{Im}\ z)^{s} \ge e^{-s d(z,h)}.$$
Since $E^{l}(z,s) =
\sum\limits_{\delta \in \Gamma_{\infty} \backslash \Gamma_{l}} 
(\operatorname{Im} \delta z )^{s}$
is a sum of positive terms over 
$\Gamma_{\infty} \backslash \Gamma$-orbits of $z$,
we obtain
$$E^{l}(z,s) \ge e^{-s \hat{d} (z,h)},$$
where $\hat{d} (z,h)$ denotes the distance from $h$ to the $\Gamma$-orbits of $z$.
We should recall two facts here.
The first one is $(2.1)$. 
The second one is that the diameters of the {\it thick parts} of $S_{l}$ are bounded by 
a positive constant $D$ for all small $L_{l}(\gamma_{1}), L_{l}(\gamma_{2})$,
where $D$ depends only on the degenerating family $\{S_{l} \}.$ 
(For example, by using the Bers constant we can easily see the second fact. Refer to
Theorem~5.2.6 in \cite{Bu} p.130.)
Then, for $z \in \partial N_{\gamma_{1}}(k) \cap S_{l}^{2}$, we can observe that
$\hat{d} (z,h) \leq  -(1+k) \log L_{l}(\gamma_{1}) + D'$.
Here $D'$ is a constant depending only on the degenerating family. Then we have
$$E^{l}(z,s) \ge e^{-s \hat{d} (z,h)}
\ge e^{-sD'}l_{1}^{(1+k)s}. $$
The remaining two cases are similar.
\end{proof}

\begin{lmm}
Let $ s > 1$. 
 For $i=1,2$, 
let $\Omega_{i}$ be any region $(\Subset S_{0}^{i+1})$
containing $\partial C_{i}(1) \cap S_{0}^{i+1}$. There exist positive constants 
$P_{i}=P_{i} ( s,\Omega_{i}, \{S_{l} \})$ depending only on $s$ and $\Omega_{i}$
 and the degenerating family
$\{S_{l} \}$ such that for any sufficiently small  $L_{l}(\gamma_{1})$, then
\begin{align*}
& (f^{l})^{*} E^{l}(z,s)  \ge P_{1} \  L_{l}(\gamma_{1})^{2s} \qquad \hspace{37pt} 
 \text{on} \ \Omega_{1},\\
& (f^{l})^{*} E^{l}(z,s) \ge P_{2}\ L_{l}(\gamma_{1})^{2s}\
 L_{l}(\gamma_{1})^{2s} \quad \; \; \;
 \text{on} \ \Omega_{2}.
\end{align*}
\end{lmm}
\begin{proof}
We will show only the first case. The second case is similar.
We set
$$P_{l} = \inf\limits_{z \in \Omega_{1}} L_{l}(\gamma_{1})^{-2s}(f^{l})^{*} E^{l}(z,s) .$$
Suppose that there exists a subsequence 
 $l_{j}\rightarrow 0$ such that
 $\lim\limits_{j \rightarrow \infty }P_{l_{j}} = 0$.
Consider the function 
$P_{l_{j}}^{-1} L_{l_{j}}(\gamma_{1})^{-2s}(f^{l_{j}})^{*} E^{l_{j}}(z,s)$.
By $(2.2)$ and $(2.3)$, we can observe that
$$( (f^{l_{j}})^{*} \Delta_{l_{j}} - s(s-1))\ P_{l_{j}}^{-1} L_{l_{j}}(\gamma_{1})^{-2s}
(f^{l_{j}})^{*}E^{l_{j}}(z,s) = 0$$
and
$$\inf\limits_{z \in \Omega_{1}} P_{l_{j}}^{-1} L_{l_{j}}(\gamma_{1})^{-2s}
(f^{l_{j}})^{*} E^{l_{j}}(z,s) = 1. $$
We choose another region $\Omega_{1}'$ such that
$ \Omega_{1} \Subset \Omega_{1}' \Subset S_{0}^{2}$.
 Because $( (f^{l_{j}})^{*} \Delta_{l_{j}} - s(s-1))$ are uniformly non-degenerate on $\Omega_{1}'$, the Harnack inequality provides
\begin{align*}
\sup\limits_{z \in \Omega_{1}} P_{l_{j}}^{-1} L_{l_{j}}(\gamma_{1})^{-2s}
(f^{l_{j}})^{*} E^{l_{j}}(z,s) 
& \leq c(\Omega_{1},\Omega_{1}')
\inf\limits_{z \in \Omega_{1}} P_{l_{j}}^{-1} L_{l_{j}}(\gamma_{1})^{-2s}
(f^{l_{j}})^{*} E^{l_{j}}(z,s)\\
&= c(\Omega_{1},\Omega_{1}') < \infty.
\end{align*}
Then using the interior Schauder estimate and the diagonal method as in the proof of
Theorem 1 in \cite{O},
we can have a further subsequence which will be denoted by the same notation such that
$P_{l_{j}}^{-1} L_{l_{j}}(\gamma_{1})^{-2s} (f^{l_{j}})^{*} E^{l_{j}}(z,s)$ and its 
first and second derivatives converge uniformly 
on any compact subset of $\Omega_{1}$ to a nonnegative function $G(z,s)$ and its derivatives
respectively.
Then $G(z,s)$ satisfies
$$(\Delta_{0}-s(s-1))\ G(z,s) = 0$$
and
$$\sup\limits_{z \in \Omega_{1}} G(z,s) \leq
\varlimsup\limits_{j \rightarrow \infty } \sup\limits_{z \in \Omega_{1}}
P_{l_{j}}^{-1} L_{l_{j}}(\gamma_{1})^{-2s} (f^{l_{j}})^{*} E^{l_{j}}(z,s) 
\leq c(\Omega_{1},\Omega_{1}')
< \infty. $$

Now it should be noted that 
$\Omega_{1} \supset (f^{l})^{-1} (\partial N_{\gamma_{1}}(1) \cap 
S_{l}^{2})$ for any sufficiently small $l$ because 
$(f^{l})^{-1} (\partial N_{\gamma_{1}}(1) \cap 
S_{l}^{2})$ converges to  $\partial C_{1} (1) \cap S_{0}^{2}$ as $l  \rightarrow 0.$
We choose another region $\Omega_{1}'' \Subset \Omega_{1}$
such that $\Omega_{1}'' \supset (f^{l})^{-1} (\partial N_{\gamma_{1}}(1) \cap 
S_{l}^{2})$ for any sufficiently small $l$.
Then we have
\begin{align*}
\sup\limits_{z \in \Omega_{1}} G(z,s) 
&\geq \sup\limits_{z \in \Omega_{1}''} G(z,s) \\
& = \lim\limits_{j \rightarrow \infty } \sup\limits_{z \in \Omega_{1}''}
P_{l_{j}}^{-1} L_{l_{j}}(\gamma_{1})^{-2s} (f^{l_{j}})^{*} E^{l_{j}}(z,s) \\
& \geq \lim\limits_{j \rightarrow \infty } \sup\limits_{z \in 
(f^{l})^{-1} (\partial N_{\gamma_{1}}(1) \cap 
S_{l}^{2})}
P_{l_{j}}^{-1} L_{l_{j}}(\gamma_{1})^{-2s} (f^{l_{j}})^{*} E^{l_{j}}(z,s) \\
& \geq \lim\limits_{j \rightarrow \infty } \inf\limits_{z \in 
(f^{l})^{-1} (\partial N_{\gamma_{1}}(1) \cap 
S_{l}^{2})}
P_{l_{j}}^{-1} L_{l_{j}}(\gamma_{1})^{-2s} (f^{l_{j}})^{*} E^{l_{j}}(z,s)\\
& = \lim\limits_{j \rightarrow \infty } \inf\limits_{w \in 
\partial N_{\gamma_{1}}(1) \cap S_{l}^{2}}
P_{l_{j}}^{-1} L_{l_{j}}(\gamma_{1})^{-2s} E^{l_{j}}(w,s)\\
& \geq \lim\limits_{j \rightarrow \infty } P_{l_{j}}^{-1} K_{1}
= +\infty. 
\end{align*}
Here we used the first inequality in Lemma 2.6.
This is a contradiction.
\end{proof}

\begin{rem}
All claims in Lemmas 2.6, 2.7 remain valid even in the case where 
the components $S_{0}^{2},S_{0}^{3}$ have other {\it old} punctures.
However, care for such additional {\it old} punctures is needed in the proof of the following theorem.
\end{rem}

\begin{thm}
Let $s > 1$. 
For all sufficiently small $ L_{l}(\gamma_{1}), L_{l}(\gamma_{2})$, it holds that
\begin{align*} 
&  E^{l}(z,s)  \ge Q_{1} \  L_{l}(\gamma_{1})^{2s} \qquad \hspace{42pt} 
 \text{on} \ S_{l}^{2} - f^{l}(C_{2}(a)),\\
 &  E^{l}(z,s)  \ge Q_{2} \ L_{l}(\gamma_{1})^{2s} L_{l}(\gamma_{2})^{s} \qquad \hspace{10pt} 
 \text{on}\  S_{l}^{2} \cap N_{\gamma_{2}}(1),\\
&  E^{l}(z,s) \ge Q_{3}\ L_{l}(\gamma_{1})^{2s}  L_{l}(\gamma_{2})^{2s}  \; \; \quad \quad  
\text{on} \ S_{l}^{3}.
\end{align*}
Here $Q_{1}=Q_{1}(s,a,\{S_{l} \})$ is a positive constant depending only on
$s, a$ and the degenerating family $\{S_{l} \}$.
$Q_{i}= Q_{i}(s, \{S_{l} \})\; (i=2,3)$ are positive constants depending only on 
$s$ and the degenerating family $\{S_{l} \}$.
\end{thm}
\begin{rem}
In the case where $S_{l}^{i}$ has additional {\it old} punctures $(i=2,3)$,
we have to replace $S_{l}^{i}$ with 
$S_{l}^{i} - f^{l}(\text{the union of all neighborhoods of old punctures}),$
and all $Q_{i}$'s depend on all the removed neighborhoods.
\end{rem}
\begin{proof}
First, we will show the first inequality. Set
$$Q_{l} = \inf\limits_{z \in (f^{l})^{-1}(N_{\gamma_{1}}(1)) \cap S_{0}^{2}} L_{l}(\gamma_{1})^{-2s}
(f^{l})^{*} E^{l}(z,s) 
= \inf\limits_{w \in N_{\gamma_{1}}(1) \cap \bar{S}_{l}^{2}} L_{l}(\gamma_{1})^{-2s}E^{l}(w,s) .$$
Due to the first inequality in Lemma 2.7, all we have to prove is that $Q_{l}$ is larger 
than a positive constant for all small $l$.
Assume that there exists a subsequence 
$l_{j} \rightarrow 0$ such that
$\lim\limits_{j \rightarrow \infty }Q_{l_{j}} = 0$.
For each $j$, we can find a point $w_{j} \in N_{\gamma_{1}}(1) \cap \bar{S}_{l}^{2}$ such that 
$$L_{l_{j}}(\gamma_{1})^{-2s}E^{l_{j}}(w_{j},s) =
\inf\limits_{w \in N_{\gamma_{1}}(1) \cap \bar{S}_{l_{j}}^{2}} L_{l_{j}}(\gamma_{1})^{-2s}E^{l_{j}}(w,s).$$
If $w_{j}\; \text{is not on the geodesic}\ \gamma_{1}$, 
set $z_{j} = (f^{l_{j}})^{-1}(w_{j})$.
Divide our situation into three cases.
(If necessary, we will take a subsequence which is denoted by the same symbol, for simplicity.)

I. infinitely many $w_{j}$ are on the geodesic $\gamma_{1}$,

II. there exists $b \ge 1$ such that all but finitely many $z_{j}$ 
are outside of $C_{1}(b) \cap S_{0}^{2},$

III. there exists a subsequence such that 
$\lim\limits_{j \rightarrow \infty} z_{j} = q_{1}.$

In case I, due to the first inequality with $k=0$ in Lemma 2.6,
$$Q_{l_{j}} = L_{l_{j}}(\gamma_{1})^{-2s}E^{l_{j}}(w_{j},s) 
\ge K_{1} L_{l_{j}}(\gamma_{1})^{-s} \ge K_{1} > 0 \quad  \text{for all large }j.$$
This is a contradiction.

In case II, we choose a region 
$\Omega_{1}$ ($\Subset  S_{0}^{2}$) 
which contains $\partial C_{1}(1) \cap S_{0}^{2}$
and $z_{j}$.
Due to the first inequality in Lemma 2.7,
$$Q_{l_{j}} = L_{l_{j}}(\gamma_{1})^{-2s} (f^{l_{j}})^{*} E^{l_{j}}(z_{j},s)  
\ge 
\inf\limits_{z \in \Omega_{1}} L_{l_{j}}(\gamma_{1})^{-2s}
(f^{l_{j}})^{*} E^{l_{j}}(z,s) \ge P_{1}(\Omega_{1}) > 0.$$
This is a contradiction.

In case III, there exists $0 \le k_{j} \le 1$ such that
$f^{l_{j}}(z_{j}) \in \partial N_{\gamma_{1}}(k_{j}) \cap S^{2}_{l_{j}}.$
Then due to the first inequality in Lemma 2.6, we have

\begin{align*} 
Q_{l_{j}} = L_{l_{j}}(\gamma_{1})^{-2s} (f^{l_{j}})^{*} E^{l_{j}}(z_{j},s) 
& \ge
 \inf\limits_{w \in \partial N_{\gamma_{1}}(k_{j}) \cap S_{l_{j}}^{2}} L_{l_{j}}(\gamma_{1})^{-2s}
 E^{l_{j}}(w,s)\\
 & \ge K_{1} L_{l_{j}}(\gamma_{1})^{(k_{j}-1)s} \ge  K_{1} > 0
 \end{align*} 
for all large $j$. This is a contradiction.
We have proved the first inequality.

Next, we will show the second inequality in a similar method.
We set
$$Q'_{l} 
= \inf\limits_{w \in N_{\gamma_{2}}(1) \cap \bar{S}_{l}^{2}} L_{l}(\gamma_{1})^{-2s} L_{l}(\gamma_{2})^{-s}
E^{l}(w,s) .$$
Assume that there exists a subsequence 
$l_{j} \rightarrow 0$ such that
$\lim\limits_{j \rightarrow \infty }Q'_{l_{j}} = 0$.
For each $j$, we can find a point $w_{j} \in N_{\gamma_{2}}(1) \cap \bar{S}_{l}^{2}$ such that 
$$ L_{l_{j}}(\gamma_{1})^{-2s}  L_{l_{j}}(\gamma_{2})^{-s} E^{l_{j}}(w_{j},s) =
\inf\limits_{w \in N_{\gamma_{2}}(1) \cap \bar{S}_{l_{j}}^{2}} L_{l_{j}}(\gamma_{1})^{-2s}
L_{l_{j}}(\gamma_{2})^{-s} E^{l_{j}}(w,s).$$
If $w_{j}\; \text{is not on the geodesic}\ \gamma_{2}$, 
set $z_{j} = (f^{l_{j}})^{-1}(w_{j})$.
Divide our situation into three cases.
(If necessary, we will take a subsequence which is denoted by the same symbol, for simplicity.)

$\operatorname{I}'.$ infinitely many $w_{j}$ are on the geodesic $\gamma_{2}$,

$\operatorname{II}'.$ there exists $b \ge 1$ such that all but finitely many $z_{j}$ 
are outside of $C_{2}(b) \cap S_{0}^{2}$,

$\operatorname{III}'.$ there exists a subsequence such that 
$\lim\limits_{j \rightarrow \infty} z_{j} = q_{2}$.

In case $\operatorname{I}'$, due to the second inequality with $k=0$ in Lemma 2.6,
$$Q'_{l_{j}} = L_{l_{j}}(\gamma_{1})^{-2s} L_{l_{j}}(\gamma_{2})^{-s} E^{l_{j}}(w_{j},s) 
\ge K_{2} > 0 \quad  \text{for all large }j.$$
This is a contradiction.

In case $\operatorname{II}'$, we choose a region 
$\Omega'_{1}$ ($\Subset S_{0}^{2}$) 
which contains $\partial C_{1}(1) \cap S_{0}^{2}$
and $z_{j}$.
Due to the first inequality in Lemma 2.7,
\begin{align*}
Q'_{l_{j}} = L_{l_{j}}(\gamma_{1})^{-2s} L_{l_{j}}(\gamma_{2})^{-s} (f^{l_{j}})^{*} E^{l_{j}}(z_{j},s)  
& \ge 
\inf\limits_{z \in \Omega'_{1}} L_{l_{j}}(\gamma_{1})^{-2s} L_{l_{j}}(\gamma_{2})^{-s}
(f^{l_{j}})^{*} E^{l_{j}}(z,s) \\
& \ge P_{1}(\Omega'_{1}) L_{l_{j}}(\gamma_{2})^{-s}
\ge P_{1}(\Omega'_{1}) > 0
\end{align*}
for all large $j$.
This is a contradiction.

In case $\operatorname{III}'$, there exists $0 \le k_{j} \le 1$ such that
$f^{l_{j}}(z_{j}) \in \partial N_{\gamma_{2}}(k_{j}) \cap S^{2}_{l_{j}}.$
Then due to the second inequality in Lemma 2.6, we have
\begin{align*}
Q'_{l_{j}} = L_{l_{j}}(\gamma_{1})^{-2s} L_{l_{j}}(\gamma_{2})^{-s} E^{l_{j}}(f^{l_{j}}(z_{j}),s) 
& \ge
\inf\limits_{w \in \partial N_{\gamma_{2}}(k_{j}) \cap S_{l_{j}}^{2}} L_{l_{j}}(\gamma_{1})^{-2s}
 L_{l_{j}}(\gamma_{2})^{-s} E^{l_{j}}(w,s)\\
& \ge K_{2} L_{l_{j}}(\gamma_{1})^{-k_{j} s} \ge  K_{2} >0
\end{align*}
for all large $j$. This is a contradiction.
We have proved the second inequality.
We can prove the third inequality in the same way as the first inequality,
using the third inequality in Lemma 2.6 and the second inequality in Lemma 2.7.
\end{proof}


\end{large}
\end{document}